\documentclass{amsart}
\usepackage{pstricks}

\usepackage{graphicx}
\newcommand{\Z}{\mathbb {Z}}

 \newtheorem{theorem}{Theorem}[section]

\theoremstyle{definition}

\theoremstyle{remark}

\numberwithin{equation}{section}

\title[Extended periodic Links and HOMFLYPT polynomial]{Extended periodic Links and HOMFLYPT polynomial}
\author{Nafaa Chbili and Hajer Jebali}
\address{Department of Mathematical Sciences, United Arab Emirates University, 15551 Al Ain, UAE}
\email{nafaachbili@uaeu.ac.ae}
\address{D\'epartement de Math\'ematiques, Facult\'e des Sciences, Universi\'e de Monastir, Monastir 5000, Tunisia}
\email{hajer.jebali@fsm.rnu.tn}

\begin{document}

\thanks{N. Chbili, Supported by an Individual Research Grant from the College of Science, United Arab Emirates University}

\maketitle

\begin{abstract}
Extended strongly periodic links have been introduced by Przytycki and Sokolov as a   symmetric surgery presentation  of three-manifolds on which the finite cyclic group acts without fixed points.
The purpose of this paper is to prove that the symmetry of these links is reflected by the first coefficients of the  HOMFLYPT  polynomial.
{\bf AMS Key words.} Periodic links, three-manifolds, HOMFLYPT polynomial.
\end{abstract}

\section{Introduction}

Let $G$ be the cyclic group of order $p$, where $p$ is a prime and let $S^3$ be the three-dimensional sphere.
We say that $G$ acts semi-freely on $S^3$, if the action of $G$ leaves a circle $\gamma$ as the set of fixed points.
By the positive solution of the Smith conjecture \cite{BM}, we know that $\gamma$ is unknotted and
such action is conjugate to an orthogonal action, which is actually  a rotation by an angle of   $\frac{2\pi}{p}$. Moreover,  the quotient space in that case is homeomorphic to $S^3$.
A link $L$ in $S^3$ is said to be periodic if it is set-wise invariant by the rotation and $L\cap \gamma =\emptyset $.  If $L$ is a periodic  link,  then we denote by $\overline L$ its  quotient (factor) link.
 For each component $\bar l_i$ of $\bar L$, let $k_i$ be the linking number between $\overline  {l_i}$ and $\gamma$. It is well known that the total  number of
 components of $L$ is $\sum_{i} gcd(k_i,p)$. A $p$-periodic link is said to be \emph{ strongly periodic} if  the linking number of each component of the quotient link with $\gamma$ is congruent to zero modulo $p$.
 The importance of this class of links is due to the fact  that they provide an equivariant surgery presentation for branched cyclic covers  along knots. More precisely,
Przytycki and Sokolov \cite{PS} proved that a three-manifold admits a semi-free action of the finite cyclic group of order $p$ with a circle as the set of fixed points if and only if $M$ is obtained from the
 three-sphere by surgery along a strongly $p-$periodic link. The study of the quantum invariants of strongly periodic links  led to many interesting   criteria   for  periodicity of three-manifolds. For instance,   the computation of the second coefficient of the Conway polynomial of strongly periodic links was used to prove that the Casson-Walker-Lescop invariants of periodic 3-manifolds satisfies a certain congruence relation \cite{Ch3}. The study of the Jones polynomial of this class of links was used to find necessary conditions for a 3-manifold to be periodic, see \cite{Ch1,Ch2,CL,BGP} for instance.\\
In the same context, Prytycki and Sokolov introduced a symmetric surgery presentation for freely-periodic 3-manifolds. More precisely, they proved that if a  3-manifold admits an action of the cyclic group  $G$ which
has no fixed points, then $M$ is obtained from $S^3$ by surgery along a link $L'=L \cup \gamma$ where $L$ is a strongly $p$-periodic  link. We will refer to the link $L'$ as  \emph{extended strongly $p-$periodic link}.
 This paper initiates the  study of the quantum  invariants of  this class of links.  Our    ultimate goal is to establish new obstruction criteria for  3-manifolds periodicity based on the study of quantum invariants of extended strongly periodic links.\\
The HOMFLYPT polynomial  is a topological invariant of isotopy classes of oriented  links that was introduced as a natural generalization of the Jones polynomial, \cite{HOMFLY,PT}. It is a  two-variable polynomial $P_L(v,z)$ that can be defined  recursively  as follows:
$$ \begin{array}{ll}
{\bf (i)}&P_{\bigcirc}(v,z)=1, \\
{\bf (ii)}&v^{-1}P_{L_{+}}(v,z)-v P_{L_{-}}(v,z) =zP_{L_{0}}(v,z),
\end{array}$$
where  $ \large {\bigcirc}$ is the trivial knot, $L_{+}$,
$L_{-}$ and  $L_{0}$ are three link projections which are identical except at a given crossing where they look as pictured below.\\

\begin{center}
\includegraphics[width=9cm,height=2.5cm]{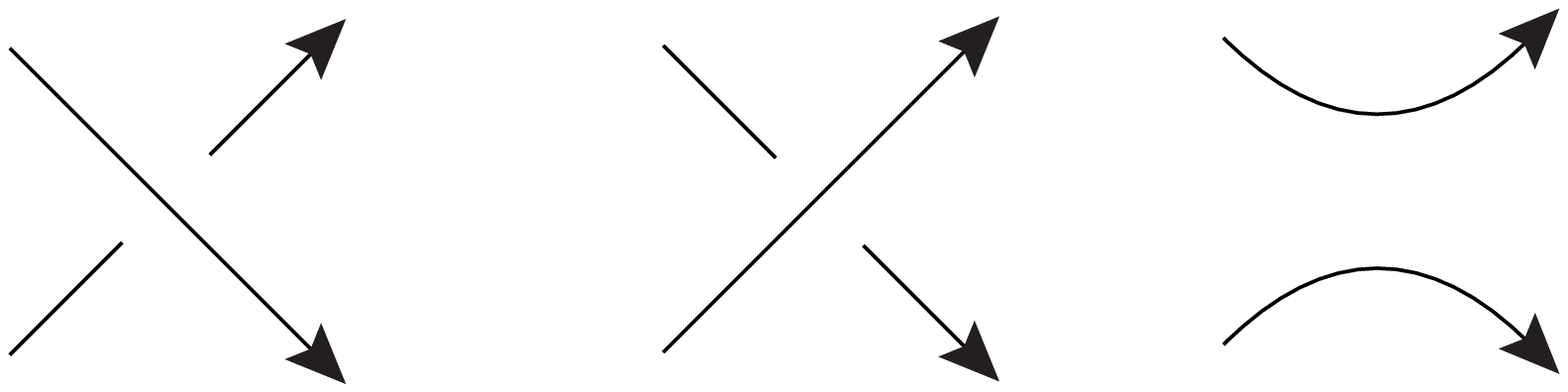}
\end{center}
$\hspace{80pt} L_{+}\hspace{90pt} L_{-}\hspace{70pt} L_{0}$
\begin{center} {\sc Figure 1.}
\end{center}

If $L$ is an $n$-components  link, then $P_L(v,z)=\displaystyle\sum_{i\geq 0} P_{1-n+2i}(v)z^{1-n+2i}$, where $P_{1-n+2i}(v) \in \Z[v^{\pm 1}]$. In particular, if $L$ is a knot then $P_L(v,z)=\displaystyle\sum_{i\geq 0} P_{2i}(v)z^{2i}$, where $P_{2i}(v) \in \Z[v^{\pm 2}]$. The purpose of this paper is to prove the following necessary conditions for a link to be an extended strongly periodic link: \\

\begin{theorem} Let $p$ be an odd prime. Let $L'=L\cup \gamma$ be an extended  strongly $p-$periodic link. If $L'$ has $n$ components,  then
\begin{enumerate}
\item $P_{1-n,L'}(v) \in \Z_{p}[v^{\pm p}]$,
\item $P_{1-n+2i,L'}(v) \equiv 0, \mbox{ mod } p$, for all $1\leq i \leq \frac{p-1}{2}. $
\end{enumerate}
\end{theorem}



Using Knotinfo \cite{CL}, we have checked that only 14 links  with a number of crossings less than or equal to 11 satisfy Condition 2 of Theorem 1.1, for $p=3$. These links are:\\
$L10n100\{0,1,0\}; L10n100\{0,1,1\}; L10n103\{0,0,0\}$ and $L10n103\{0,0,1\}$, which have
$P_{-1}(v)=-3v^{-7}+6v^{-5}-3v^{-3}$,\\
$L10n107\{0,0,0\}; L10n107\{1,0,0\} ; L10n107\{0,1,0\}; L10n107\{1,1,0\};L10n107\{0,0,1\},\\
L10n107\{1,0,1\};L10n107\{0,1,1\}$ and $L10n107\{1,1,1\}$ which  have $P_{-1}(v)=0$,
$L11n448\{1,1,0\}$ and $L11n448\{1,1,1\}$ which have  $P_{-1}(v)=3v^3-6v^5+3v^7$.\\
 Since our ultimate goal is the study of the   quantum invariants of freely periodic three-manifolds,  we would like to mention that one can write a colored version  of Theorem 1.1. In other words,  these results can be
   be generalized to the link $L'^r=L\cup \gamma^{r}$ which is obtained by taking a diagram of $L'$, then  substituting  $\gamma$  by  $r$ parallel copies. The consequences of this colored version on the quantum invariants of periodic 3-manifolds  will be discussed  in a forthcoming paper.\\
   
\begin{theorem} Let $p$ be an odd prime. If $L'^r$ has $n$ components,  then
\begin{enumerate}
\item $(v-v^{-1})^{1-r}P_{1-n,L'^r}(v) \in \Z_{p}[v^{\pm p}]$,
\item $P_{1-n+2i,L'^r}(v) \equiv 0, \mbox{ mod } p$, for all $1\leq i \leq \frac{p-1}{2}. $
\end{enumerate}
\end{theorem}
\section{Periodic links}
 In this section, we shall  define   extended periodic  links  and  briefly review  some of their properties which will be needed in rest of the paper.  \\

\textbf{Definition 2.1.} {\sl
 Let $p\geq 2$ be an integer. A link $L$ of $S^{3}$ is said to be
$p$-periodic if and only if there exists an orientation-preserving
auto-diffeomorphism $h$ of $S^{3}$
such that:\\
\begin{tabular}{rl}
&1) $\gamma=$Fix($h$) is homeomorphic to the  circle $S^{1}$,\\
&2) the link  L is disjoint from $\gamma$,\\
&3) $h$ is of order $p$,\\
&4) $h(L)=L$.
\end{tabular}\\
  If $L$ is  periodic we will denote the quotient link by
$\overline L$.\\}

By the positive solution of the Smith Conjecture \cite{BM}, $\gamma$ is unknotted and the  diffeomorphism $h$ is conjugate to a rotation  $\varphi_p$ by an angle of $\frac{2\pi}{p}$.
Recall here that if the quotient link $\overline L$ is a knot then
the link $L$ may have more than one component. In general, the
number of components of $L$ depends on the linking numbers of the
components of $\overline L$ with the axis
of the rotation.\\
Note that if $L$  is a $p$-periodic link, then there exists an  $n$-tangle $T$, such that $L$ is isotopic to the closure of $T^p$;   $L=\widehat{T^p}$. Consequently, $L$ has a diagram which is invariant by a planar rotation of angle $\frac{2\pi}{p}$.
 A special class of periodic links has been introduced by Przytycki and Sokolov \cite{PS} as in the following definition.

{\bf Definition 2.2.} {\sl  Let $p\geq 2$ be an integer. A
$p$-periodic link $L$ is said to be strongly $p$-periodic if and
only if one of the following equivalent conditions holds:\\
\begin{tabular}{rl}
&(i) The linking number of each component of $\overline L$ with
the axis $\Delta$ is 0 modulo $p$.\\
 &(ii) The group $\Z_p$ acts freely on the set of components of
 $L$.\\
 &(iii) The number of components of $L$ is $p$ times greater than
 the number of\\
 & components  of  $
\overline L$.
\end{tabular}
} \\

If $L$ is a strongly $p$ periodic link, then the link $L'=L\cup \gamma$ is called an extended strongly $p$-periodic link.  Notice that
 an extended  strongly $p$-periodic link $L'$ has $p\alpha+1$ components, where $\alpha$ is the
 number of components of the quotient link $\overline L$. The
 $p\alpha$ components are divided into $\alpha$ orbits with
 respect to the free action of  $\Z_p$ (condition (i)).\\
 Assume that  $\overline L=l_1\cup \dots \cup\l_{\alpha}$, there is a natural
cyclic order on each orbit of components of $L$. Namely,
$$L=l_1^1 \cup\dots \cup l_1^p \cup l_2^1 \cup\dots \cup l_2^p \cup \dots \cup l_{\alpha}^1 \cup\dots
\cup l_{\alpha}^p$$ where $ \varphi_p(l_i^t)=l_{i}^{t+1} \;
\forall 1\leq t\leq p-1$ and $\varphi_p(l_i^p)=l_{i}^{1}$, for all
$1\leq i \leq \alpha$.\\

\begin{center}
\includegraphics[width=4cm,height=2.5cm]{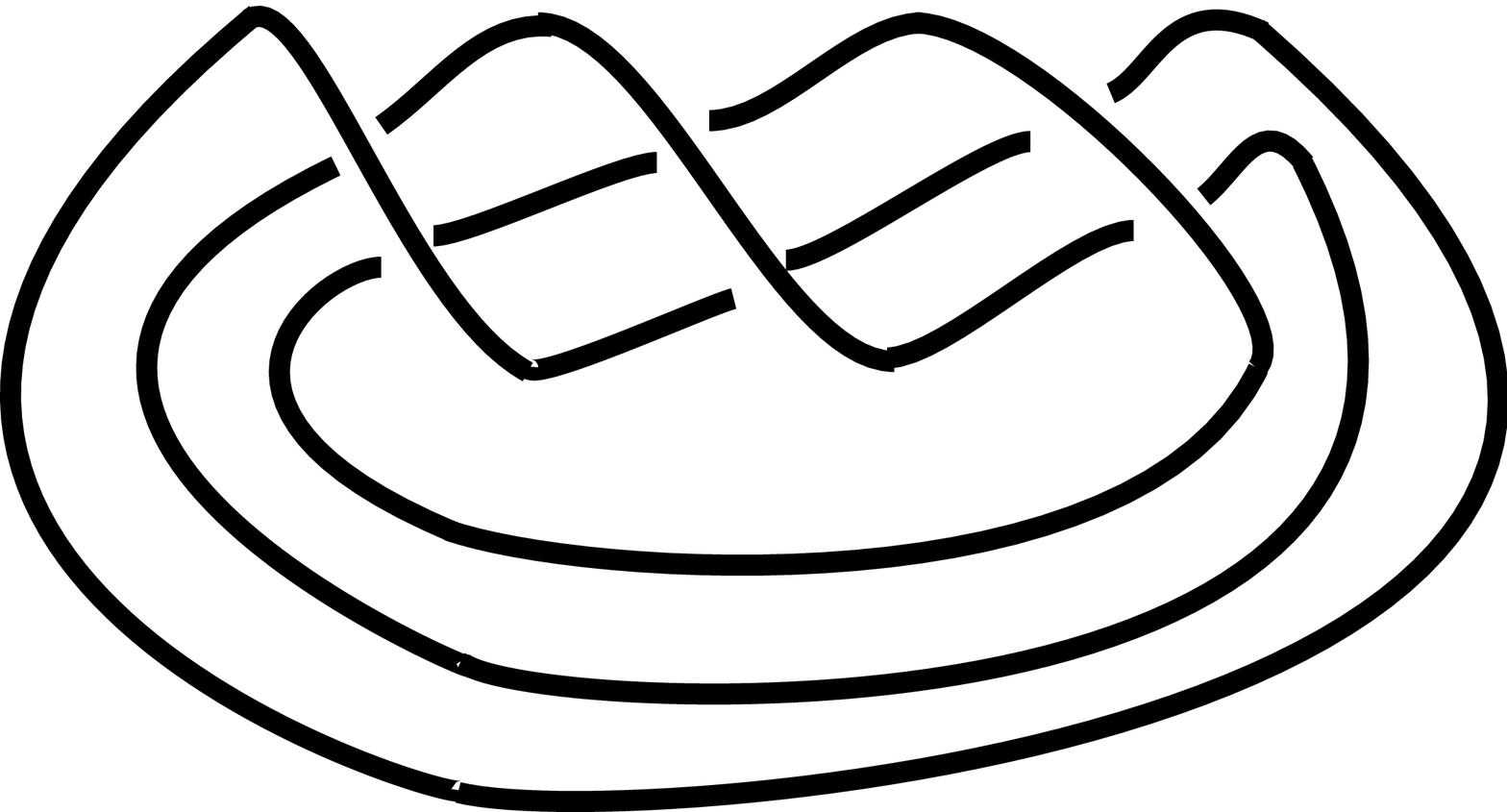} \hspace{70pt}
\includegraphics[width=4cm,height=2.5cm]{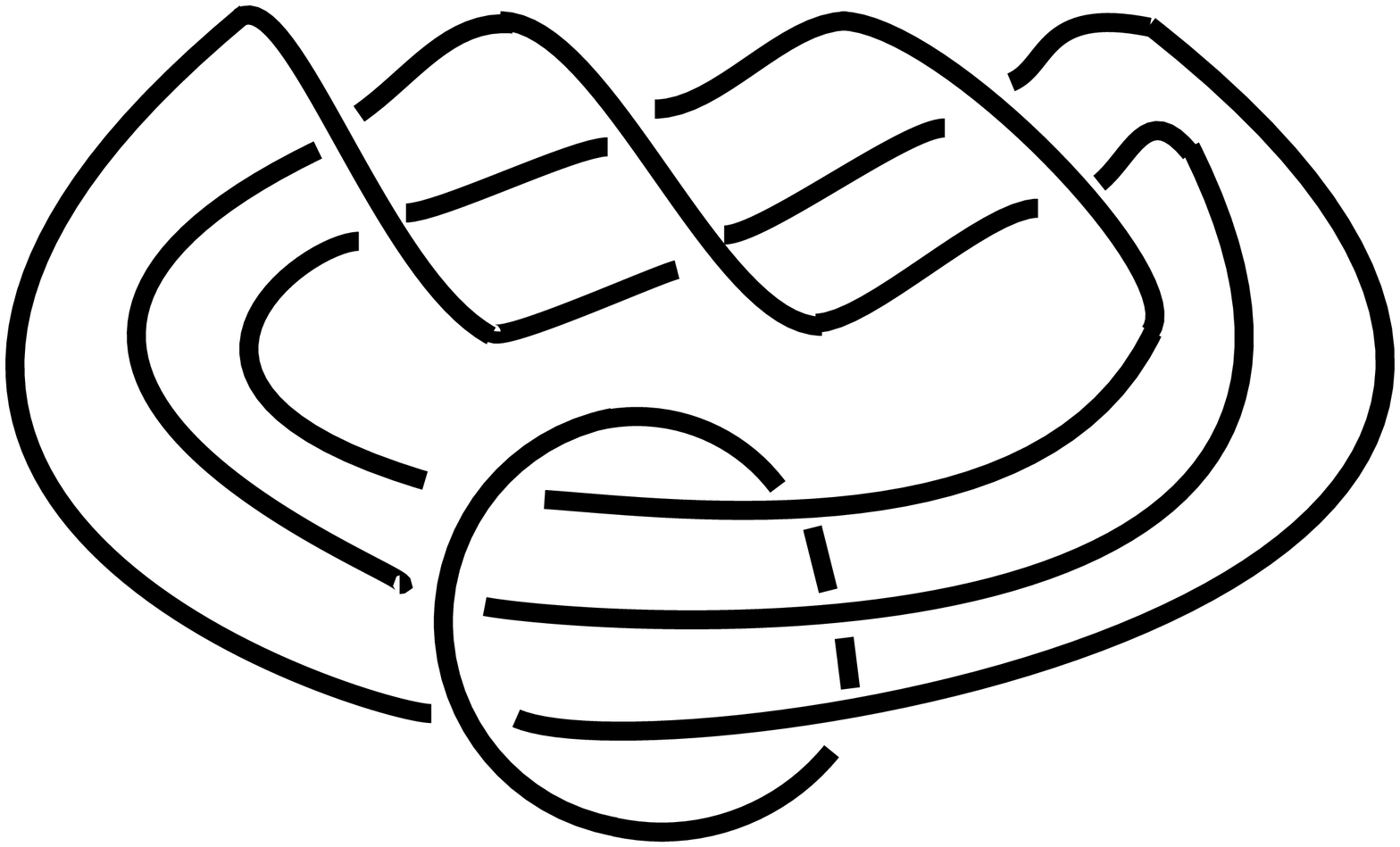}\\
{\sc Figure 2:} A strongly 3-periodic link in left and an extended strongly 3-periodic link in the right
\end{center}

 A three-manifold $M$ is said to be $p$-periodic
 if the finite cyclic group  $\Z_p$ acts semi-freely on $M$
 with a circle as the set of fixed points. According to  \cite{PS} and \cite{Sa1}, we know
  that a three-manifold is $p$-periodic if and only if is
  obtained from $S^{3}$ by surgery along a
strongly $p$-periodic link.  A three-manifold is said to be freely $p$-periodic  if it admits a fixed-points free  action of    $\Z_p$. Przytycki and Sokolov  gave an
equivariant surgery presentation of such manifolds. More precisely, they  proved that a three-manifold is freely $p$-periodic if and only if it can be obtained from $S^{3}$ by surgery along an extended
strongly $p$-periodic link \cite{PS}.

\section{Proofs}

\textbf{Proof of Theorem 1.1}\\

The proof of  Condition 1 is straightforward using the following Lemma: \\

{\bf Lemma 3.1  \cite{LM}.} {\sl Let $K=k_{1}\cup
k_{2}\cup\dots\cup k_{n}$ be an $n$-component link. Then
$$P_{1-n,K}(v)=v^{2\lambda}(v^{-1}-v)^{n-1}\displaystyle\prod_{i=1}^{n}P_{0,k_{i}}(v),$$
where $\lambda$  denotes the total linking number of the link
$K$.}

In the case of an extended strongly $p$-periodic link, $L'=l_1^1 \cup\dots \cup l_1^p \cup l_2^1 \cup\dots \cup l_2^p \cup \dots \cup l_{\alpha}^1 \cup\dots
\cup l_{\alpha}^p \cup \gamma$, we get:

$$P_{-p\alpha,L'}(v)=v^{2\lambda}(v^{-1}-v)^{p\alpha}P_{0,\gamma}(v)\displaystyle\prod_{i=1}^{\alpha}\prod_{k=1}^{p}P_{0,l_{i}^k}(v)=v^{2\lambda}(v^{-1}-v)^{p\alpha}\displaystyle\prod_{i=1}^{\alpha}[P_{0,l_{i}^1}(v)]^p$$
We conclude using the fact that the total linking number $\lambda$ is a multiple of $p$.\\

Now, we shall prove Condition 2 of  Theorem 1.1.  It is worth mentioning that the proof is a generalization of  some  techniques   used in \cite{Tr1,Tr2,Yo} in the case of periodic knots.
 We will first introduce some notation.
Let $L$ be a  $p$-periodic link in the three-sphere  and let   $\overline
L$ be its quotient link. Thus  $L=\pi^{-1}(\overline L)$,
where $\pi$ is the canonical surjection corresponding to the
action of the rotation $\varphi_p$ on the three-sphere. Let $\overline L_+$,
$\overline L_-$ and $\overline L_0$ denote the three links which
are identical to $\overline L$ except near one crossing where they
are like in Figure 1. Now, let $L_{p+}:=\pi^{-1}(\overline L_+)$,
$L_{p-}:=\pi^{-1}(\overline L_-)$ and $L_{p0}:=\pi^{-1}(\overline
L_0)$.  Let $D_{p+}$ (respectively $D_{p-}$ and $D_{p0}$) a  link diagram
 of $L_{p+}$ (respectively $L_{p-}$, $L_{p0}$) which is invariant by a planar rotation.
 We define an equivariant crossing change as a change that transforms
  $D_{p+}$ into $D_{p-}$ or vice-versa. If this crossing change
 involves two different components of the factor  link, then we call it a
mixed equivariant crossing change. Otherwise, it is called a self
equivariant crossing change.\\

The proof of Condition 2 will be done by induction on the number of crossings of $D$. If $D$ has no crossings, then the condition is obviously true as  $D'$ is  the unlink with  $p\alpha+1$ components.\\
Now,  our strategy is as follows. Starting from  a diagram of an extended strongly periodic link $D'=D \cup \gamma$,  we shall   use equivariant crossing changes to
reduce $D'$ to an extended torus link $T(p,kp) \cup \gamma$  keeping the coefficients $P_{1-n+2i}(v)$ unchanged modulo $p$ and up to the multiplication by  a power of $v^{2p}$, for  $1\leq i \leq \frac{p-1}{2}.$
Without loss of generality, we can assume that $D'=D_{p+}'$.  First, we prove the following lemma.\\

\textbf{Lemma 3.2.} $v^{-p}P_{D'_{p+}}(v,z)-v^{p}P_{D'_{p-}}(v,z)\equiv z^pP_{D'_{p0}}(v,z) \mbox{\;\; modulo } p.$\\

\emph{Proof.}  The idea of the proof is based on  writing  the skein relations along $p$ crossings of $L$ which belong to the same orbit.
 We notice that since $p$ is prime then  the contribution of certain links will add to zero  modulo $p$. We remain only with the extended strongly periodic links
   $D'_{p+}$, $D'_{p-}$ and $D'_{p0}$ and  it can be seen that their polynomials satisfy the congruence in the lemma, see also \cite{Tr1}.\\

 \textbf{Lemma 3.3.}  { \sl For  all  $0\leq i \leq \frac{p-1}{2}$, we have:   $P_{-p\alpha+2i, D'_{p+}}(v) \equiv v^{2p}P_{-p\alpha + 2i, D'_{p-}}(v) $  modulo $p$.}\\

 \emph{Proof.} Assume that $D'_{p+}$ has $p\alpha+1$ components, then $D'_{p-}$
 has $p\alpha+1$ components as well.  To determine  the number of components of  link  $D'_{p0}$  we shall  look at the effect of the crossing change at  the factor link.
 If the crossing to be smoothed  is mixed, then two components will merge into one  whose linking number with $\gamma$ is zero modulo $p$, as the sum of the two linking numbers.
  Hence,  $D'_{p0}$  will be an extended strongly $p$-periodic link with $p(\alpha-1)+1$ components.  Consequently,  the lowest power of $z$ that appears in the right hand side of the congruence in Lemma 3.2 is
  $z^{p}z^{-p(\alpha-1)}=z^{-p\alpha+2p}$. Hence,
 $v^{-p}P_{-p\alpha +2i, D'_{p+}}(v) -v^{p} P_{-p\alpha +2i, D'_{p-}}(v) \equiv 0 $ modulo $p$ for $0 \leq i \leq p-1 $ and Lemma 3.3 is proved in this case.\\

 If the crossing to be smoothed in the factor link  is a self-crossing, then that component will split into two components and  we will have two cases.\\
 Case 1:  $D'_{p0}$  is  an extended periodic link with $p(\alpha-1)+3$.  In this case the lowest power of $z$ that appears in the right hand side of the congruence in Lemma 2 is $z^{p}z^{1-(p(\alpha-1)+3)}=z^{-p\alpha+2p-2}$.
 Consequently, \\
  $v^{-p}P_{-p\alpha +2i, D'_{p+}}(v) -v^{p} P_{-p\alpha +2i, D'_{p-}}(v) \equiv 0 $ modulo $p$ for $0 \leq i < p-1$.\\
 Case 2: $D'_{p0}$ is an extended  strongly $p$-periodic link with $p(\alpha+1)+1$ components. In this case, the diagram $D'_{p0}$ has less crossings than $D'_{p+}$. Then we can use the  induction assumption to conclude.

By the previous Lemma one can use the equivariant crossing change  without changing the first coefficients of the HOMFLYPT polynomial modulo $p$ and up to the multiplication by a power of $v^p$.
 If the factor diagram $\overline{L}=l_1 \cup l_2  \dots \cup  l_{\alpha}$ has  more than one component,
then we can  change some mixed  crossings in order to put the component   $l_1$ over the other  components. Repeating this process for the remaining components, we should be able to split the factor link diagram into disjoint components.
 The corresponding equivariant crossing changes will then split  the diagram of $D_{p+}$ into a disjoint  union of $\alpha$ orbits.  The diagram of $D'_{p+}$ can be then seen as a connected sum of  extended strongly periodic link diagrams each of which has exactly $p+1$ components. Consequently, it will be enough for the rest of the proof to  assume that $D'_{p+}$  has exactly $p+1$ components. \\

\textbf{Lemma 3.4.} {\sl The diagram  $D'_{p+}=D_{p+}\cup \gamma$ can be transformed into a trivial  $(p+1)$-component link or  $T(p,kp) \cup \gamma $ by equivariant crossing changes.}\\
\emph{Proof.} The diagram $D_{p+}$ can be seen as the closure of an $n$-tangle of the form  $T^p$, where the closure of $T$ represents the factor link diagram $\overline L$.
Using equivariant crossing changes,  $D_{p+}$ can be transformed into the closure of a braid  of the form $B^p$ or a trivial $p$-component link, see \cite{Tr1,Ch3}.
Following \cite{Tr1}, we can use equivariant crossing changes  to  transform  the closure of $B^p$ into a torus link of the type  $T(p,kp)$. Then, $D'_{p+}$ is transformed into $T(p,kp) \cup \gamma$.\\

\textbf{Lemma 3.5.}  If  $L'=T(p,pk)\cup \gamma$, then $P_{-p+2i,L'}(v) \equiv 0, \mbox{ mod } p$, for all $1\leq i \leq \frac{p-1}{2}. $  \\

\emph{Proof.}  Recall that $L'$ is made up of $p$ components permuted cyclically by the action of the rotation, in addition to the unknotted component $\gamma$.
  Consequently, if we make a crossing change between $\gamma$ and any of the other
 components then the resulting links will be isotopic. Thus, if we write the skein relation corresponding to the change of these $p$ crossings, we get:\\
 $v^{-p}P_{L',p+}(v,z)-v^{p}P_{L',p-}(v,z) \equiv z^pP_{L',p0}(v,z)$. Since $L'_{p0}$ is a knot, we can see that
  $v^{-p}P_{L',p+}(v,z)-v^{p}P_{L',p-}(v,z) \equiv 0$ modulo $p$ for $1\leq 2i <p$.\\
   We can repeat these changes if necessary until we  split $\gamma$ from the torus link. Then $P_{L'}(v,z)$ is congruent to $P_{T(p,kp)}(v,z)$  modulo $p$ and up to the multiplication by  a power of $v^{2p}$.
   The latter satisfies Condition 2 of Theorem 1.1 as it was proved in   Lemma 2.8 of \cite{Yo}.\\

\textbf{Proof of Theorem 1.2}\\
The condition 1 can be proved by using
 Lemma 3.1 for an extended strongly $p$-periodic link, $L'^r=l_1^1\cup\cdots\cup l_1^p\cup l_2^1\cup\cdots\cup l_2^p\cup\cdots\cup l_\alpha^1\cup\cdots\cup l_\alpha^p\cup\gamma^r$. We have:
\begin{align*}
P_{1-p\alpha-r,L'^r}(v)&=\displaystyle v^{2\lambda}(v^{-1}-v)^{p\alpha+r-1}\left(P_{0,\gamma}(v)\right)^r\prod_{i=1}^{\alpha}\left(P_0,l_i^{1}(v)\right)^p\\
&=v^{2\lambda}(v^{-1}-v)^{p\alpha+r-1}\prod_{i=1}^{\alpha}\left(P_{0,l_i^1}(v)\right)^p\,.
\end{align*}

Now it is sufficient to use the fact that the total linking number $\lambda$ is a multiple of $p$.\\

Let us  prove   condition 2 of Theorem 1.2.  Let $D'^r=D\cup \gamma^r$ be a diagram of the  extended strongly periodic link $L'^r$.
The proof is done by induction on $r$. The case $r=1$ corresponds to Theorem 1.1. Assume  that $r >1$ and that the condition holds for $r-1$.  By the equivariant crossing changes illustrated above, we should be able to split the orbits of $D$ as it was done in the proof of Lemma 3.3. Then, transform each of these orbits into either a trivial $p$-component link or a torus link of type $T(p,kp)$ as explained in the proof of Lemma 3.4. In the next step we can use crossing changes between one of the copies of   $\gamma$ and the torus link to split them from each other. Repeating these operations if necessary, we can split that copy of $\gamma$ from the rest of link. Then, we can conclude using the induction assumption.

\end{document}